\documentclass[12pt,twoside]{article}
\usepackage{amsmath,amsfonts,amssymb}
\usepackage{graphicx}
\usepackage{eucal}
\usepackage{mathrsfs}
\usepackage{theorem}
\usepackage{pifont}
\usepackage {epsfig}
\usepackage {graphicx}
\newcommand{\Real}{\mathbb{R}}

\newcommand{\Complex}{\mathbb{C}}

\newcommand{\todo}[1]{{\sffamily To do:}}

\newtheorem{theorem}{Theorem}
\newtheorem{proposition}{Proposition}

\newtheorem {corollary}{Corollary}
\newtheorem {lemma}{Lemma}
\newenvironment{proof}{{\flushleft \emph{Proof}:}}{\ding{110}}

\title{Eigenvalue inequalities in terms of Schatten norm bounds on differences of semigroups, and application to
Schr\"odinger operators}
\author{Michael Demuth\footnotemark[1],\; Guy Katriel{\footnotemark[1]\;\;\footnotemark[2]}}
\date{}

\begin{document}

\maketitle
\renewcommand{\thefootnote}{\fnsymbol{footnote}}
\footnotetext[1]{Institute of Mathematics, Technical University of
Clausthal, 38678 Clausthal-Zellerfeld, Germany.}
\footnotetext[2]{Partially supported by the Minerva Foundation
(Germany).}
\begin{abstract}
We develop a new method for obtaining bounds on the negative
eigenvalues of self-adjoint operators $B$ in terms of a Schatten
norm of the difference of the semigroups generated by $A$ and $B$,
where $A$ is an operator with non-negative spectrum. Our method is
based on the application of the Jensen identity of complex function
theory to a suitably constructed holomorphic function, whose zeros
are in one-to-one correspondence with the negative eigenvalues of
$B$. Applying our abstract results, together with bounds on Schatten
norms of semigroup differences obtained by Demuth and Van Casteren,
to Schr\"odinger operators, we obtain inequalities on moments of the
sequence of negative eigenvalues, which are different from the
Lieb-Thirring inequalities.
\end{abstract}
%


\section{Introduction}

Let $A$ be a self-adjoint operator on a complex Hilbert space, whose
spectrum is non-negative. If $B$ is another self-adjoint operator,
such that the difference $D_t=e^{-tB}-e^{-tA}$ of the semigroups
corresponding to $A,B$ belongs to a Schatten ideal (trace class or
Hilbert-Schmidt class), we will prove inequalities which provide
bounds from above on the negative eigenvalues of $B$, in terms of
Schatten norms of $D_t$. The usefulness of such results follows from
the fact that for concrete operators, for example when $B$ is a
Schr\"odinger operator $B=-\Delta+V$, and $A$ is the free
Schr\"odinger operator $A=-\Delta$, it is known that, under
appropriate conditions on the potential $V$, $D_t$ belongs to a
Schatten ideal, and explicit bounds on the Schatten norm of $D_t$
are available \cite{demuth}. Indeed such results are important in
the study of the absolutely continuous spectrum of the perturbed
operator $B$. The theorems proven here show that these bounds on the
Schatten norms of $D_t$ can also be used in the study of the
discrete spectrum of $B$.

The method used to prove our results is based on constructing a
holomorphic function whose zeros are in one-to-one correspondence
with the negative eigenvalues of $B$, and using complex function
theory to bound these zeros. Specifically we will use the Jensen
identity (see, e.g., \cite{rudin}, p. 307):
\begin{lemma}
\label{jensen0} Let $\Omega_r$ be an open disk centered at $0$ and
with radius $r$.  Let $h:U \rightarrow \Complex$ be a holomorphic
function on the open set $U$, where ${\overline{\Omega}}_r\subset
U$, and assume $h(0)=1$. Then
$$\frac{1}{2\pi}\int_0^{2\pi}\log(|h(r e^{i\theta})|)d\theta=\log\Big(\prod_{z\in \bar{\Omega}_r, h(z)=0}
\frac{r}{|z|}\Big)=\int_0^r \frac{n(u)}{u}du,$$ where $n(u)$ ($0\leq
u\leq r$) denotes the number of zeros of $h$ in $\bar{\Omega}_u$.
\end{lemma}

In Section \ref{general} we prove general theorems which give bounds
on the moments (sums of powers) of the sequence of negative
eigenvalues of an operator $B$ in terms of the trace norm of the
semigroup difference. In Section \ref{generalhs} we prove analogous
bounds in terms of the Hilbert-Schmidt norm of the semigroup
difference. In Section \ref{schrodinger} we apply the theorems of
Section \ref{generalhs} to derive inequalities for the negative
eigenvalues of Schr\"odinger operators under some conditions on the
potential, which are different from the well-known Lieb-Thirring
inequalities.

\section{Eigenvalue inequalities in terms of trace-norm bounds on semigroup differences}
\label{general}

In this section we will prove results under the assumption that
$A$,$B$ are selfadjoint operators, with the spectrum of $A$
non-negative, and such that the difference of semigroups
$D_t=e^{-tB}-e^{-tA}$ is of trace class. This implies that the
negative spectrum of $B$, which we denote by
$$\sigma^-(B)=\sigma(B)\cap (-\infty,0),$$
 consists only of
eigenvalues, which can accumulate only at $0$ (of course compactness
of $D_t$ is sufficient for this property). We shall denote by
$N(-s)$ the number of eigenvalues $\lambda$ of $B$ which satisfy
$\lambda< -s$.

We begin by proving identities expressing the moments of the
negative eigenvalues of the operator $B$ in terms of an integral.
It should be noted that the identities hold also in the
case that one side is infinite - which implies that the other side
is infinite too.

\begin{theorem}\label{identity}
Let $A$,$B$ be self-adjoint in a complex Hilbert space ${\cal{H}}$,
with $\sigma(A)\subset[0,\infty)$. Assume that $D=e^{-B}-e^{-A}$ is
of trace class. Then, for any $\gamma>1$, we have
\begin{eqnarray}\label{idnz}
&&\sum_{\lambda\in \sigma^-(B)} |\lambda|^{\gamma}\\ &=&\frac{\gamma
(\gamma-1)}{2\pi}\int_0^1 \frac{1}{r}|\log(r)
|^{\gamma-2}\int_0^{2\pi}\log\Big(|Det(I-F(re^{i\theta}))|\Big)d\theta
dr,\nonumber
\end{eqnarray}
where $F(z)$ is the operator-valued function defined by
\begin{equation}\label{df}F(z)=z[I-ze^{-A}]^{-1}D,\end{equation}
and for $\gamma=1$ we have
\begin{equation}\label{identity1}
\sum_{\lambda\in \sigma^-(B)} |\lambda|=\lim_{r\rightarrow
1}\frac{1}{2\pi}\int_0^{2\pi}\log\Big(|Det(I-F(re^{i\theta}))|\Big)d\theta.
\end{equation}
\end{theorem}

\begin{proof}
We have, for all $z\in \Complex$
\begin{equation}\label{sta}I-ze^{-B}=I-ze^{-A}-zD,\end{equation} and if also $|z|<1$, so that
$\|ze^{-A}\|<1$ then $I-ze^{-A}$ is invertible, so that $F(z)$ given
by (\ref{df}) is well defined, and we can write (\ref{sta}) as
$$[I-ze^{-A}]^{-1}[I-ze^{-B}]=I-F(z).$$
Thus we have the following equivalence for $|z|<1$:
\begin{equation*}  \log(z)\in
\sigma(B) \;\Leftrightarrow  \frac{1}{z}\in \sigma(e^{-B})
\;\Leftrightarrow 1\in \sigma(F(z)),
\end{equation*}
so that
\begin{equation}\label{es}\sigma^-(B)=\{ \;\log(z)
\;|\;|z|<1,\;\;1\in \sigma(F(z))\}. \end{equation} Since we assume
$D$ is of trace class, then so is $F(z)$. We note also that
\begin{equation}\label{f0}
F(0)=0.
\end{equation}
Since $F(z)$ is a trace class operator, the determinant
$$h(z)=Det(I-F(z))$$
is well defined, and we have that $h$ is holomorphic in the unit
disk and
$$h(z)=0 \;\;\Leftrightarrow \;\;1\in \sigma(F(z))
\;\;\Leftrightarrow\;\;\log(z)\in \sigma^-(B).$$ Thus
$$\sigma^-(B)=\{\; \log(z)\;|\; |z|<1,\;\; h(z)=0\;\},$$
so that, for all $s>0$,
\begin{equation}\label{nn}
N(-s)=n(e^{-s}),
\end{equation}
where $n(r)$ denotes the number of zeros of $h$ in
$\Omega_r=\{z\;|\;|z|<r\}$. By (\ref{f0}) we have
$$h(0)=Det(I-F(0))=Det(I)=1.$$
Applying the Jensen identity, Lemma \ref{jensen0}, we have, for any
$0<r<1$,
\begin{equation}\label{je}\frac{1}{2\pi}\int_0^{2\pi}\log(|h(r
e^{i\theta})|)d\theta=\int_0^r \frac{n(u)}{u}du,\end{equation} and
making the substitution $u=e^{-s}$ in the integral on the
right-hand side of (\ref{je}) and using (\ref{nn}) we get
\begin{equation}\label{je2}\frac{1}{2\pi}\int_0^{2\pi}\log(|h(r
e^{i\theta})|)d\theta =\int_{\log(\frac{1}{r})}^\infty N(-s)
ds.\end{equation} We now recall the well-known identity
\begin{eqnarray}\label{wk}
\sum_{\lambda\in \sigma^-(B)} |\lambda|^{\gamma}=\gamma\int_0^\infty
s^{\gamma-1}N(-s)ds.
\end{eqnarray}
Taking $\gamma=1$, (\ref{wk}) becomes
\begin{eqnarray}\label{wk1}
\sum_{\lambda\in \sigma^-(B)} |\lambda|=\int_0^\infty N(-s)ds.
\end{eqnarray}
Taking the limit $r\rightarrow 1$ in (\ref{je2}), we have
\begin{equation}\label{je3} \lim_{r\rightarrow 1}\frac{1}{2\pi}\int_0^{2\pi}\log(|h(r
e^{i\theta})|)d\theta =\int_{0}^\infty N(-s) ds.\end{equation} From
(\ref{wk1}) and (\ref{je3}), we conclude
\begin{equation*}\label{identity01}
\sum_{\lambda\in \sigma^-(B)} |\lambda|=\lim_{r\rightarrow
1}\frac{1}{2\pi}\int_0^{2\pi}\log(|h(r e^{i\theta})|)d\theta,
\end{equation*}
so that we have (\ref{identity1}).

We now assume that $\gamma>1$. Multiplying (\ref{je2}) by
$\frac{1}{r}|\log(r)|^{\gamma-2}$ and integrating over $r\in [0,1]$,
we obtain
\begin{eqnarray*}
&&\frac{1}{2\pi}\int_0^1\int_0^{2\pi}\frac{1}{r}|\log(r)|^{\gamma-2}\log(|h(r
e^{i\theta})|)d\theta dr\nonumber\\
&=&\int_0^1
\frac{1}{r}|\log(r)|^{\gamma-2}\int_{\log(\frac{1}{r})}^\infty N(-s)
ds dr\nonumber\\&=&\int_0^\infty  N(-s)\int_{e^{-s}}^1
\frac{1}{r}|\log(r)|^{\gamma-2} dr
ds=\frac{1}{\gamma-1}\int_0^\infty N(-s)s^{\gamma-1}
ds,\end{eqnarray*} which, together with (\ref{wk}), implies
\begin{eqnarray*}
\sum_{\lambda\in \sigma^-(B)}|\lambda|^{\gamma}=
\frac{\gamma(\gamma-1)}{2\pi}\int_0^1\int_0^{2\pi}\frac{1}{r}
|\log(r)|^{\gamma-2}\log(|h(r e^{i\theta})|)d\theta dr,\nonumber
\end{eqnarray*}
so that we have (\ref{idnz}).
\end{proof}

By bounding the function $h$ of Theorem \ref{identity} from above,
we obtain bounds on the moments of the negative eigenvalues.

\begin{theorem}\label{ggiq}
Let $A$,$B$ be self-adjoint in a complex Hilbert space ${\cal{H}}$,
with $\sigma(A)\subset[0,\infty)$. Assume that $D=e^{-B}-e^{-A}$ is
of trace class. Then for any $\gamma>1$,
\begin{equation}\label{giniq}
\sum_{\lambda\in \sigma^-(B)} |\lambda|^{\gamma}\leq \frac{\gamma
(\gamma-1)}{2\pi}\int_0^1 |\log(r)|^{\gamma-2}\int_0^{2\pi}
\|[I-re^{i\theta}e^{-A}]^{-1}D \|_{tr}d\theta dr,
\end{equation}
and for $\gamma=1$ we have
\begin{equation*}
\sum_{\lambda\in \sigma^-(B)} |\lambda|\leq\limsup_{r\rightarrow
1}\frac{1}{2\pi}\int_0^{2\pi}\|[I-re^{i\theta}e^{-A}]^{-1}D\|_{tr}d\theta.
\end{equation*}
\end{theorem}

\begin{proof}
We recall the general inequality for trace class operators $T$ (see,
e.g., \cite{simonb})
\begin{equation}\label{indet}|Det(I-T)|\leq e^{\|T\|_{tr}},\end{equation}
which gives
\begin{eqnarray*}
\log\Big(|Det(I-F(re^{i\theta}))|\Big)\leq \|F(re^{i\theta})\|_{tr}=
r\|[I-re^{i\theta}e^{-A}]^{-1}D\|_{tr}.
\end{eqnarray*}
Substituting this inequality into (\ref{idnz}), (\ref{identity1}),
we obtain the results.
\end{proof}

Bounding the integral on the right-hand side of (\ref{giniq}), we
get the following theorem. Although we shall later prove a stronger
result, Theorem \ref{exp}, it is useful to present Theorem
\ref{prim}, whose proof is more straightforward, and for which the
coefficient in the inequalities can be evaluated explicitly, in
terms of Euler's $\Gamma$-function and Riemann's $\zeta$-function.

\begin{theorem}\label{prim}
Let $A$,$B$ be self-adjoint in a complex Hilbert space ${\cal{H}}$,
with $\sigma(A)\subset[0,\infty)$. Assume that, for some $t>0$,
$D_t=e^{-tB}-e^{-tA}$ is of trace class.

Then, for any $\gamma> 2$, we have the inequality
\begin{equation}\label{pini}
\sum_{\lambda\in \sigma^-(B)} |\lambda|^{\gamma}\leq
\Gamma(\gamma+1)\zeta(\gamma-1)\frac{1}{t^\gamma}\|D_t\|_{tr},
\end{equation}
and the right-hand side is finite.
\end{theorem}

\begin{proof}
We note first that it suffices to prove (\ref{pini}) for $t=1$, that
is, setting $D=D_1=e^{-B}-e^{-A}$, to prove
\begin{equation}\label{pinis}
\sum_{\lambda\in \sigma^-(B)} |\lambda|^{\gamma}\leq
\Gamma(\gamma+1)\zeta(\gamma-1)\|D\|_{tr},
\end{equation}
since (\ref{pini}) follows from (\ref{pinis}) by replacing $A,B$ by
$tA,tB$.

Since $\sigma(A)\subset [0,\infty)$, we have $\|e^{-A}\|\leq 1$, so
that, for $|z|<1$,
\begin{eqnarray*}
\|[I-ze^{-A}]^{-1}\|\leq \frac{1}{1-|z|},
\end{eqnarray*}
hence
\begin{eqnarray*}
\|F(r e^{i\theta})\|_{tr}=r\|[I-r
e^{i\theta}e^{-A}]^{-1}D\|_{tr}&\leq& r\|[I-r
e^{i\theta}e^{-A}]^{-1}\| \|D\|_{tr}\nonumber\\ &\leq& \|D\|_{tr}
\frac{r}{1-r}.
\end{eqnarray*}
From the inequality (\ref{giniq}) of Theorem \ref{ggiq} we thus have
\begin{eqnarray}\label{a0}
&&\sum_{\lambda\in \sigma^-(B)}|\lambda|^{\gamma}\\&\leq&
\frac{\gamma (\gamma-1)}{2\pi}\int_0^1
|\log(r)|^{\gamma-2}\int_0^{2\pi} \|[I-re^{i\theta}e^{-A}]^{-1}D
\|_{tr}d\theta dr\nonumber
\\&\leq& \gamma
(\gamma-1)\|D\|_{tr}\int_0^1 |\log(r)|^{\gamma-2} \frac{1}{1-r}
dr\nonumber
\\&=& \gamma
(\gamma-1)\|D\|_{tr}\int_0^{\infty}\frac{x^{\gamma-2}}{e^{x}-1}
dx=\Gamma(\gamma+1)\zeta(\gamma-1)\|D\|_{tr}\nonumber
\end{eqnarray}
so we have (\ref{pinis}).
\end{proof}

A more refined estimate on the integral in (\ref{giniq}) yields the
following theorem, which is stronger than Theorem \ref{prim}. We
note that this theorem is valid for $\gamma>1$, rather than
$\gamma>2$ as in Theorem \ref{prim}. The value of the constant
$C_{tr}(\gamma)$ is given in the proof of the theorem, in terms of
some integrals.

\begin{theorem}\label{exp}
Let $A$,$B$ be self-adjoint in a complex Hilbert space ${\cal{H}}$,
with $\sigma(A)\subset[0,\infty)$. Assume that, for some $t>0$,
$D_t=e^{-tB}-e^{-tA}$ is of trace class.

Then, for any $\gamma> 1$, we have the inequality
\begin{equation}\label{ini}
\sum_{\lambda\in \sigma^-(B)} |\lambda|^{\gamma}\leq
C_{tr}(\gamma)\frac{1}{t^\gamma}\|D_t\|_{tr},
\end{equation}
where $C_{tr}(\gamma)$ is a finite constant depending only on
$\gamma$.
\end{theorem}

\begin{proof}
As noted in the proof of Theorem \ref{prim}, it suffices to prove
(\ref{ini}) for $t=1$, that is, setting $D=e^{-B}-e^{-A}$, to prove
\begin{equation}\label{inis}
\sum_{\lambda\in \sigma^-(B)} |\lambda|^{\gamma}\leq
C_{tr}(\gamma)\|D\|_{tr}.
\end{equation}

Since $\sigma(e^{-A})\subset [0,1]$, we have
\begin{eqnarray}\label{bi2}
&&|z|\|[I-ze^{-A}]^{-1}\|=\|[z^{-1}I-e^{-A}]^{-1}\|\nonumber\\&\leq
&\frac{1}{\min_{u\in [0,1]}|z^{-1}-u|}= \left\{\begin{array}{cc}
                                              \frac{1}{|z^{-1}-1|}, &  Re(z^{-1})\geq 1\\
                                              \frac{1}{|Im(z^{-1})|}, & 0< Re(z^{-1})< 1\\
                                              \frac{1}{|z^{-1}|}, & Re(z^{-1})\leq 0
                                              \end{array}\right .
\end{eqnarray}
so that
\begin{eqnarray*}
\|F(re^{i\theta})\|_{tr}&=&r\|[I-re^{i\theta}e^{-A}]^{-1}D\|_{tr}\leq
r\|[I-re^{i\theta}e^{-A}]^{-1}\| \|D\|_{tr} \nonumber\\&\leq&
r\|D\|_{tr}\left\{
\begin{array}{cc}
                                              \frac{1}{\sqrt{r^2-2r\cos(\theta)+1}} & \cos(\theta)\geq r\\
                                              \frac{1}{|\sin(\theta)|} & 0<\cos(\theta)<r\\
                                              1 & \cos(\theta)\leq 0
                                              \end{array}\right.
\end{eqnarray*}
From the inequality (\ref{giniq}) of Theorem \ref{ggiq} we thus
have, for $\gamma>1$,
\begin{eqnarray}\label{a00}
&&\sum_{\lambda\in \sigma^{-}(B)}|\lambda|^{\gamma}\\&\leq&
\frac{\gamma (\gamma-1)}{2\pi}\int_0^1
|\log(r)|^{\gamma-2}\int_0^{2\pi} \|[I-re^{i\theta}e^{-A}]^{-1}D
\|_{tr}d\theta dr\nonumber
\\&\leq&\frac{\gamma
(\gamma-1)}{\pi}\|D\|_{tr}\Big[
\int_{0}^{1}|\log(r)|^{\gamma-2}\int_0^{\arccos(r)}
\frac{1}{\sqrt{r^2-2r\cos(\theta)+1}} d\theta dr\nonumber
\\&+&\int_{0}^1
|\log(r)|^{\gamma-2}\int_{\arccos(r)}^{\frac{\pi}{2}} \frac{1
}{|\sin(\theta)|} d\theta dr+ \int_0^1 \int_{\frac{\pi}{2}}^\pi
|\log(r)|^{\gamma-2} d\theta dr\Big]\nonumber.
\end{eqnarray}
We estimate the integrals in (\ref{a00}) from above: making the
substitution $s=\frac{1}{r}$, $y=\frac{1}{\cos(\theta)}$, we have
\begin{eqnarray*}c_1(\gamma)&=&\int_{0}^{1}|\log(r)|^{\gamma-2}\int_0^{\arccos(r)}
\frac{1}{\sqrt{r^2-2r\cos(\theta)+1}} d\theta dr\nonumber
\\
&=&\int_{1}^{\infty}
\frac{(\log(s))^{\gamma-2}}{s\sqrt{s^2+1}}\int_{1}^{s}
\frac{1}{\sqrt{y-\frac{2s}{s^2+1}}}\frac{1}{\sqrt{y^2-1}}\frac{1}{\sqrt{y}} dyds\nonumber\\
&\leq&\int_{1}^{\infty}
\frac{(\log(s))^{\gamma-2}}{s\sqrt{s^2+1}}\int_{1}^{s}
\frac{1}{\sqrt{y-\frac{2s}{s^2+1}}}\frac{1}{\sqrt{y-1}} dyds\nonumber\\
&=&\int_{1}^{\infty}
\frac{(\log(s))^{\gamma-2}}{s\sqrt{s^2+1}}\log\Big(\frac{(\sqrt{s(s+1)}+\sqrt{s^2+1})^2}{s-1}
\Big)ds\\
\end{eqnarray*}
and since, for any $\epsilon>0$, the integrand in the last integral
is $O((s-1)^{\gamma-2-\epsilon})$ as $s\rightarrow 1+$ and
$O(s^{\epsilon-2})$ as $s\rightarrow \infty$, this integral is
finite whenever $\gamma>1$, so $c_1(\gamma)$ is finite.
\begin{eqnarray*}c_2(\gamma)&=&\int_{0}^1
|\log(r)|^{\gamma-2}\int_{\arccos(r)}^{\frac{\pi}{2}}
\frac{1 }{|\sin(\theta)|} d\theta dr\\
&=&\frac{1}{2}\int_{0}^1
|\log(r)|^{\gamma-2}\log\Big(\frac{1+r}{1-r}\Big)dr,
\end{eqnarray*}
and since, for any $\epsilon>0$, the integrand in the last integral
is $O(r^{1-\epsilon})$ as $r\rightarrow 0$ and
$O((1-r)^{\gamma-2-\epsilon})$ as $r\rightarrow 1$, this integral is
finite whenever $\gamma>1$, so $c_2(\gamma)$ is finite. Finally, we
have
$$c_3(\gamma)= \int_0^1 \int_{\frac{\pi}{2}}^\pi
|\log(r)|^{\gamma-2} d\theta dr=\frac{\pi}{2}\int_0^\infty
x^{\gamma-2}e^{-x} dx=\frac{\pi}{2}\Gamma(\gamma-1).$$ From
(\ref{a00}) we thus obtain, for $\gamma>1$,
\begin{equation*}
\sum_{\lambda\in \sigma^-(B)}|\lambda|^{\gamma}\leq
\frac{1}{\pi}\gamma(\gamma-1)[c_1(\gamma)+c_2(\gamma)+c_3(\gamma)]\|D\|_{tr}.
\end{equation*}
so that we have (\ref{inis}), with
\begin{equation*}C_{tr}(\gamma)=\frac{1}{\pi}\gamma(\gamma-1)[c_1(\gamma)+c_2(\gamma)+c_3(\gamma)].
\end{equation*}
\end{proof}
\newpage
One could ask what is the {\it{best}} constant $C_{tr}(\gamma)$ in
inequality (\ref{ini}), that is, given $\gamma>1$, what is the
smallest number $C_{tr}(\gamma)$ for which (\ref{ini}) will hold for
{\it{any}} pair of selfadjoint operators with $\sigma(A)\subset
[0,\infty)$. We do not know how to answer this question, but we can
give a simple {\it{lower bound}} for the possible values of
$C_{tr}(\gamma)$. We recall that the Lambert W-function is defined
on  $[-e^{-1},\infty)$  as the inverse of the function $f(x)=xe^x$.
\begin{proposition}
If $\gamma>1$ and $C_{tr}(\gamma)$ is a constant for which Theorem
\ref{exp} holds, then
\begin{equation}\label{lb}C_{tr}(\gamma)\geq -W(-\gamma
e^{-\gamma})(\gamma+W(-\gamma
e^{-\gamma}))^{\gamma-1}.\end{equation}
\end{proposition}
\begin{proof}
If the inequality (\ref{ini}) holds then in particular it must hold
when $A$ and $B$ are $1\times 1$ matrices. Thus let $A=0$, $B=-b$,
($b> 0$), $t=1$. Then the moment of the negative eigenvalues of $B$
of order $\gamma$ is simply $b^\gamma$, and
$D_1=e^{-B}-e^{-A}=e^{b}-1$. Thus inequality (\ref{ini}) becomes in
this case
$$b^\gamma\leq C_{tr}(\gamma)(e^{b}-1).$$
Since this must hold for all $b>0$, we have
\begin{equation}\label{we}C_{tr}(\gamma)\geq
\sup_{b>0}\frac{b^\gamma}{e^{b}-1}.\end{equation} By differentiating
the function of $b$ on the right-hand side of (\ref{we}) we find its
maximum on $[0,\infty)$ to be given by the expression on the
right-hand of (\ref{lb}).
\end{proof}

As an example, we take $\gamma=2$. From (\ref{lb}) we obtain
$C_{tr}(2)\geq 0.647..$ Theorem \ref{exp} gives (evaluating the
integrals numerically) $C_{tr}(2)\leq 2.5..$

Using the above argument one can see that for $\gamma<1$, Theorem
\ref{exp} {\it{cannot}} be true. Indeed, if $\gamma<1$, then the
expression on the right-hand side of (\ref{we}) goes to $+\infty$ as
$b\rightarrow 0$, so that the supremum is infinite.

We remark that the inequalities for the moments of eigenvalues
derived here imply inequalities for the {\it{number}} of eigenvalues
less than a given negative number $-s$ ($s>0$), which we denote by
$N(-s)$. Indeed since
$$\sum_{\lambda\in \sigma^-(B)}|\lambda|^{\gamma}\geq \sum_{\lambda\in \sigma(B)\cap
(-\infty,-s)}|\lambda|^{\gamma}\geq \sum_{\lambda\in \sigma(B)\cap
(-\infty,-s)}s^{\gamma}=s^\gamma N(-s),$$ we have, from (\ref{ini}),
assuming that $D_t$ is trace-class for all $t>0$,
$$N(-s)\leq \frac{1}{s^\gamma}\sum_{\lambda\in \sigma(B)\cap
(-\infty,0)}|\lambda|^{\gamma}\leq \inf_{t>0,
\gamma>1}\frac{C_{tr}(\gamma)}{(st)^\gamma}\|D_t\|_{tr}.$$

\section{Eigenvalue inequalities in terms of Hilbert-Schmidt norm bounds on semigroup differences}
\label{generalhs}

In this section we prove theorems analogous to those in the previous
section, for the case in which the semigroup difference is
Hilbert-Schmidt rather than trace class. The proofs are similar, the
difference being that we have to get around the fact that the
determinant is not defined for a general Hilbert-Schmidt
perturbation of the identity. In the applications to Schr\"odinger
operators, it is easier to verify that the semigroup difference is
Hilbert-Schmidt than to verify that it is trace class, so the
theorems of this section will be used in these applications, to be
presented in Section \ref{schrodinger}.

The following theorem is the Hilbert-Schmidt analog of Theorem
\ref{identity}. It should be noted, however, that unlike in Theorem
\ref{identity}, here we have only inequalities rather than
identities.

\begin{theorem}\label{identityhs}
Let $A$,$B$ be self-adjoint in a complex Hilbert space ${\cal{H}}$,
with $\sigma(A)\subset[0,\infty)$. Assume that $D=e^{-B}-e^{-A}$ is
Hilbert-Schmidt. Then we have, for any $\gamma>1$,
\begin{eqnarray}\label{idnzhs}
&&\sum_{\lambda\in \sigma^-(B)} |\lambda|^{\gamma}
\\&\leq &\frac{\gamma (\gamma-1)}{2\pi}\int_0^1
\frac{1}{r}|\log(r)
|^{\gamma-2}\int_0^{2\pi}\log\Big(|Det(I-(F(re^{i\theta}))^2)|\Big)d\theta
dr.\nonumber
\end{eqnarray}
where $F(z)$ is the operator-valued function defined by
\begin{equation*}F(z)=z[I-ze^{-A}]^{-1}D,\end{equation*}
and for $\gamma=1$
\begin{equation}\label{identity1hs}
\sum_{\lambda\in \sigma^-(B)} |\lambda|\leq \lim_{r\rightarrow
1}\frac{1}{2\pi}\int_0^{2\pi}\log\Big(|Det(I-(F(re^{i\theta}))^2)|\Big)d\theta.
\end{equation}
\end{theorem}

\begin{proof}
Like in the proof of Theorem \ref{identity}, we have
$$\sigma^-(B)=\{\; \log(z) \;|\;|z|<1,\;\;1\in \sigma(F(z))\;\}.$$
Since we assume $D$ is Hilbert-Schmidt, then so is $F(z)$,  and this
implies that $(F(z))^2$ is trace class, so we can define the
holomorphic function
$$h(z)=Det(I-(F(z))^2),$$
and we have
$$1\in \sigma(F(z))\;\;\Rightarrow \;\; 1\in \sigma((F(z))^2)\;\;\Leftrightarrow\;\;h(z)=0,$$
and thus
\begin{equation}\label{equi}\sigma^-(B)\subset \{\; \log(z)\;|\;
|z|<1,\;\;h(z)=0\;\}.\end{equation} Since (\ref{equi}) is an inclusion rather
than an equality as in (\ref{es}), (\ref{nn}) is replaced by the
inequality
\begin{equation*}
N(-s)\leq n(e^{-s}),
\end{equation*}
Since $F(0)=0$ we have $h(0)=1$. Applying the Jensen identity, as in
the proof of Theorem \ref{identity}, we get the results.
\end{proof}

The next theorem is the Hilbert-Schmidt analog of Theorem
\ref{ggiq}.
\begin{theorem}\label{ineqhs}
Let $A$,$B$ be self-adjoint in a complex Hilbert space ${\cal{H}}$,
with $\sigma(A)\subset[0,\infty)$. Assume that $D=e^{-B}-e^{-A}$ is
Hilbert-Schmidt. Then, for any $\gamma> 1$, we have the inequality
\begin{equation}\label{gini1}
\sum_{\lambda\in \sigma^-(B)} |\lambda|^{\gamma}\leq \frac{\gamma
(\gamma-1)}{2\pi}\int_0^1 r|\log(r)|^{\gamma-2}\int_0^{2\pi}
\|[I-re^{i\theta}e^{-A}]^{-1}D \|_{HS}^2d\theta dr.
\end{equation}
and for $\gamma=1$ we have
\begin{equation*}
\sum_{\lambda\in \sigma^-(B)} |\lambda|\leq\limsup_{r\rightarrow
1}\frac{1}{2\pi}\int_0^{2\pi}\|[I-re^{i\theta}e^{-A}]^{-1}D\|_{HS}^2d\theta.
\end{equation*}
\end{theorem}

\begin{proof}
Using (\ref{indet}), we have
\begin{eqnarray}\label{ei010}
\log\Big(|Det(I-(F(z))^2)|\Big)\leq \|
(F(z))^2\|_{tr},
\end{eqnarray}
and since, for any Hilbert-Schmidt operator $T$ we have
$\|T^2\|_{tr}\leq \|T\|_{HS}^2$, we get
\begin{equation}\label{hhk}\|(F(z))^2\|_{tr}\leq
\|F(z)\|_{HS}^2.\end{equation} From (\ref{ei010}) and (\ref{hhk}),
together with (\ref{idnzhs}), (\ref{identity1hs}), we obtain the
results.
\end{proof}

The following theorem is the Hilbert-Schmidt analog of Theorem
\ref{exp}.
\begin{theorem}\label{exphs}
Let $A$,$B$ be self-adjoint in a complex Hilbert space ${\cal{H}}$,
with $\sigma(A)\subset[0,\infty)$. Assume that, for some $t>0$,
$D_t=e^{-tB}-e^{-tA}$ is Hilbert-Schmidt.

Then, for every $\gamma>2$, we have the inequality
\begin{equation}\label{ini1}
\sum_{\lambda \in \sigma^-(B)}|\lambda|^{\gamma}\leq
C_{HS}(\gamma)\frac{1}{t^\gamma}\|D_t\|_{HS}^2,
\end{equation}
where $C_{HS}(\gamma)$ is a finite constant depending only on
$\gamma$.
\end{theorem}

\begin{proof}
We first note that it suffices to prove (\ref{ini1}) with $t=1$,
that is, setting $D=D_1=e^{-B}-e^{-A}$, to prove
\begin{equation}\label{inig1}
\sum_{\lambda \in \sigma^-(B)}|\lambda|^{\gamma}\leq
C_{HS}(\gamma)\|D\|_{HS}^2,
\end{equation}
since (\ref{ini1}) follows from (\ref{inig1}) by replacing $A,B$ by
$tA,tB$.

Using the inequality (\ref{bi2}), we have
\begin{eqnarray*}\label{sb1}\|[I-e^{i\theta}e^{-A}]^{-1}D\|_{HS}^2 &\leq&\|[I-e^{i\theta}e^{-A}]^{-1}
\|^2 \|D\|_{HS}^2
\\&\leq& \|D\|_{HS}^2\left\{
\begin{array}{cc}
                                              \frac{1}{r^2-2r\cos(\theta)+1} & \cos(\theta)\geq r\\
                                              \frac{1}{(\sin(\theta))^2} & 0<\cos(\theta)<r\\
                                              1 & \cos(\theta)\leq 0
                                              \end{array}\right.\end{eqnarray*}
Therefore from inequality (\ref{gini1}) of Theorem \ref{ineqhs}
\begin{eqnarray*}&&\sum_{\lambda\in \sigma^-(B)
}|\lambda|^{\gamma}\\&\leq& \frac{\gamma (\gamma-1)}{2\pi}\int_0^1
r|\log(r)|^{\gamma-2}\int_0^{2\pi} \|[I-re^{i\theta}e^{-A}]^{-1}D
\|_{HS}^2d\theta dr \nonumber\\&\leq&\|D_t\|_{HS}^2\frac{\gamma
(\gamma-1)}{\pi}\Big[\int_{0}^{1}|\log(r)|^{\gamma-2}\int_0^{\arccos(r)}
\frac{r}{r^2-2r\cos(\theta)+1} d\theta dr \nonumber\\&+& \int_{0}^1
r |\log(r)|^{\gamma-2}\int_{\arccos(r)}^{\frac{\pi}{2}} \frac{1
}{(\sin(\theta))^2} d\theta dr+  \int_0^1 \int_{\frac{\pi}{2}}^\pi
r|\log(r)|^{\gamma-2} d\theta dr \Big].
\end{eqnarray*}
To verify that the above integrals are indeed finite for $\gamma>2$,
we estimate from above:
\begin{eqnarray*}c_4(\gamma)&=&\int_{0}^{1}|\log(r)|^{\gamma-2}\int_0^{\arccos(r)}
\frac{r}{r^2-2r\cos(\theta)+1} d\theta dr\\
&=&\int_{0}^{1} |\log(r)|^{\gamma-2}
\frac{2r}{1-r^2}\arctan\Big(\sqrt{\frac{1+r}{1-r}}\Big) dr,
\end{eqnarray*}
and since, for any $\epsilon>0$ the integrand is $O(r^{1-\epsilon})$
as $r\rightarrow 0$, and $O((1-r)^{\gamma-3})$ as $r\rightarrow 1$,
the integral is finite when $\gamma>2$.
\begin{eqnarray*}c_5(\gamma)&=&\int_{0}^1 r
|\log(r)|^{\gamma-2}\int_{\arccos(r)}^{\frac{\pi}{2}}
\frac{1 }{(\sin(\theta))^2} d\theta dr\\
\\&=&\int_{0}^1
|\log(r)|^{\gamma-2}\frac{r^2}{\sqrt{1-r^2}} dr,
\end{eqnarray*}
and since, for any $\epsilon>0$, the integrand is
$O(r^{2-\epsilon})$ as $r\rightarrow 0$, and
$O((1-r)^{\gamma-\frac{5}{2}})$ as $r\rightarrow 1$, the integral is
finite when $\gamma>\frac{3}{2}$. Finally,
$$c_6(\gamma)= \int_0^1 \int_{\frac{\pi}{2}}^\pi r|\log(r)|^{\gamma-2}
d\theta dr=\frac{\pi}{2}\int_0^{\infty}e^{-2x}x^{\gamma-2}dx=\pi
2^{-\gamma}\Gamma(\gamma-1),$$ finite for any $\gamma>1$. From
(\ref{gini1}) we thus have, for $\gamma>2$,
\begin{equation*}
\sum_{\lambda\in \sigma^-(B)}|\lambda|^{\gamma}\leq
\frac{1}{\pi}\gamma(\gamma-1)[c_4(\gamma)+c_5(\gamma)+c_6(\gamma)]\|D\|_{HS}^2.
\end{equation*}
so that (\ref{inig1}) holds, with
$C_{HS}(\gamma)=\frac{1}{\pi}\gamma(\gamma-1)[c_4(\gamma)+c_5(\gamma)+c_6(\gamma)].$
\end{proof}

An argument involving one-dimensional operators, like in the end of
the previous section, shows that Theorem \ref{exphs} is {\it{not}}
true if $\gamma<2$.

\section{Application to Schr\"odinger operators}
\label{schrodinger}

We now apply our general results to the study of the discrete
spectrum of Schr\"odinger operators $-\Delta+V$. Recall that the potential $V:\Real^d\rightarrow
\Real$ is said to belong to the class $K(\Real^d)$ if
$$\lim_{t\rightarrow 0}\sup_{x\in\Real^d}\int_0^{t}(e^{\eta \Delta
}|V|)(x)d\eta=0.$$ $V$ is said to belong to class $K^{loc}(\Real^d)$
if $\chi_{Q}V\in K(\Real^d)$ for any ball $Q \subset \Real^d$, where
$\chi_{Q}$ denotes the characteristic function of $Q$. $V$ is said
to be a Kato potential if $V_-=\min(V,0)\in K(\Real^d)$ and
$V_{+}=\max(V,0)\in K^{loc}(\Real^d)$.

By the min-max principle, the eigenvalues of $-\Delta+V_-$ are
smaller then or equal to the corresponding eigenvalues of
$-\Delta+V$, and therefore we have
\begin{equation}\label{pni}\sum_{\lambda \in \sigma^-(-\Delta+V)}|\lambda|^{\gamma}\leq \sum_{\lambda \in
\sigma^-(-\Delta+V_-)}|\lambda|^{\gamma},\end{equation} so that to
bound the left-hand side of (\ref{pni}) it suffices to bound the
right-hand side. We shall therefore take $A=H_0=-\Delta$,
$B=H_0+V_-$, so that
$$D_t=e^{-t(H_0+V_-)}-e^{-tH_0}.$$

We quote the following bounds for the Hilbert-Schmidt norm of $D_t$
(\cite{demuth}, Theorem 5.7)
\begin{lemma}\label{hsn} Assuming $V_-\in K(\Real^d)$, we have
\begin{equation*}\| D_{t}
\|_{HS}^2\leq 2t\int_{\Real^d}e^{-2t(H_0+V_-)}(x,x)|V_-(x)|dx
.\end{equation*}
\end{lemma}
\begin{lemma}\label{hsn22} Assuming $V_-\in K(\Real^d)$, we have
\begin{equation*}\| D_{t}
\|_{HS}^2\leq t^2\int_{\Real^d}e^{-2t
(H_0+V_-)}(x,x)|V_-(x)|^2dx.\end{equation*}
\end{lemma}
We also quote the following inequality (see \cite{demuth}, p. 66, in the proof of
Theorem 2.9):
\begin{lemma}\label{ni} Assuming $V_-\in K(\Real^d)$, we have
\begin{eqnarray*}e^{-t(H_0+V_-)}(x,y)\leq
\|e^{-t(H_0+2V_-)}\|_{L^1,L^{\infty}}^{\frac{1}{2}}(e^{-tH_0}(x,y))^{\frac{1}{2}}.
\end{eqnarray*}
\end{lemma}
Since $e^{-tH_0}(x,x)=\frac{1}{(4\pi t)^{\frac{d}{2}}},$ Lemmas
\ref{hsn},\ref{hsn22} and \ref{ni} imply
\begin{equation}\label{dv1}\| D_{t}
\|_{HS}^2\leq \frac{2t}{(8\pi
t)^{\frac{d}{4}}}\|e^{-2t(H_0+2V_-)}\|_{L^1,L^{\infty}}^{\frac{1}{2}}\|V_-\|_{L^1}
,\end{equation}
\begin{equation}\label{dv2}\| D_{t}
\|_{HS}^2\leq \frac{t^2}{(8\pi t)^{\frac{d}{4}}}
\|e^{-2t(H_0+2V_-)}\|_{L^1,L^{\infty}}^{\frac{1}{2}}
\|V_-\|_{L^2}^2.\end{equation} From (\ref{dv1}) and Theorem
\ref{exphs} we have
\begin{theorem}\label{map01}
Let $V$ be a Kato potential, and assume also $V_-\in L^1(\Real^d)$.
We have the following inequality for any $\gamma>2$,
\begin{equation*}
\sum_{\lambda \in \sigma^-(-\Delta+V)}|\lambda|^{\gamma}\leq
\frac{2C_{HS}(\gamma)}{(8\pi )^{\frac{d}{4}}}\|V_-\|_{L^1}\inf_{t>0}
\frac{\|e^{-2t(H_0+2V_-)}\|_{L^1,L^{\infty}}^{\frac{1}{2}}}{t^{\gamma+\frac{d}{4}-1}}.
\end{equation*}
\end{theorem}
Similarly, from (\ref{dv2}) and Theorem \ref{exphs} we have
\begin{theorem}\label{map02}
Let $V$ be a Kato potential, and assume also $V_-\in L^2(\Real^d)$.
We have the following inequality for any $\gamma>2$,
\begin{equation*}
\sum_{\lambda \in \sigma^-(-\Delta+V)}|\lambda|^{\gamma}\leq
\frac{C_{HS}(\gamma)}{(8\pi
)^{\frac{d}{4}}}\|V_-\|_{L^2}^2\inf_{t>0}
\frac{\|e^{-2t(H_0+2V_-)}\|_{L^1,L^{\infty}}^{\frac{1}{2}}}{t^{\gamma+\frac{d}{4}-2}}.
\end{equation*}
\end{theorem}
In order to make the bounds given by Theorems
\ref{map01},\ref{map02} more explicit we are going to bound
$\|e^{-2t(H_0+2V_-)}\|_{L^1,L^{\infty}}$ in terms of the quantity
($c>0$)
\begin{eqnarray}\label{da}\beta(c)&=&\|(c-\Delta)^{-1}V_-\|_{L^{\infty}}.\end{eqnarray}
We note that (see, {\it{e.g.}}, \cite{krishna}, Lemma 4.2.4) $V_-\in
K(\Real^d)$ implies that
\begin{equation}\label{al}
\lim_{c\rightarrow \infty}\beta(c)= 0.
\end{equation}

From \cite{demuth}, Proposition 2.2, we have
\begin{lemma}\label{kb}
Assume $V$ is a Kato potential. Then, for any $c> 0$ for which
$\beta(c)<1$, we have
$$\|e^{-t(H_0+V_-)}\|_{L^{\infty},L^{\infty}}\leq \frac{e^{ct}}{1-\beta(c)}.$$
\end{lemma}
\begin{lemma}\label{kash}Let $V$ be a Kato potential.
If $c> 0$ is such that
\begin{equation}\label{ll1}
4\beta(c)<1, \end{equation} then
$$\|e^{-2t(H_0+2V_-)}\|_{L^1,L^{\infty}}\leq \frac{1}{(4\pi
t)^{\frac{d}{2}}}\frac{e^{ct}}{1-4\beta(c)}.$$
\end{lemma}

\begin{proof}
We have (as in \cite{demuth}, proof of Theorem 2.9):
\begin{eqnarray*}&&\|e^{-2t(H_0+2V_-)}\|_{L^1,L^{\infty}}\leq
\|e^{-t(H_0+2V_-)}\|_{L^1,L^{2}}\|e^{-t(H_0+2V_-)}\|_{L^2,L^{\infty}}\nonumber\\
&=&\|e^{-t(H_0+2V_-)}\|_{L^2,L^{\infty}}^2\leq
\|e^{-t(H_0+4V_-)}\|_{L^\infty,L^{\infty}}\|e^{-tH_0}\|_{L^1,L^{\infty}}\nonumber\\
&=&\frac{1}{(4\pi
t)^{\frac{d}{2}}}\|e^{-t(H_0+4V_-)}\|_{L^\infty,L^{\infty}}.
\end{eqnarray*}
Using Lemma \ref{kb}, we get the result.
\end{proof}

Using Lemma \ref{kash}, Theorem \ref{map01} implies, for $c$
satisfying (\ref{ll1}),
\begin{equation}\label{inn11}
\sum_{\lambda \in \sigma^-(-\Delta+V)}|\lambda|^{\gamma}\leq
\frac{2^{\frac{d}{4}+1}}{(8\pi )^{\frac{d}{2}}}C_{HS}(\gamma)
\frac{1}{t^{\gamma+\frac{d}{2}-1}}\frac{e^{\frac{1}{2}ct}}{[1-4\beta(c)]^{\frac{1}{2}}}\|V_-\|_{L^1}.
\end{equation}
We can now minimize the expression on the right-hand side of
(\ref{inn11}) over $t$. Since
$$\min_{t>0}\frac{e^{\frac{1}{2}ct}}{t^{\gamma+\frac{d}{2}-1}}=\Big(\frac{e c}{2\gamma+d-2}\Big)^{\gamma+\frac{d}{2}-1}$$
we obtain
\begin{theorem}\label{mapw}
Let $V$ be a Kato potential, and assume also $V_-\in L^1(\Real^d)$.
If $c>0$ is such that $4\beta(c)<1$, then, for any $\gamma>2$,
\begin{equation*}
\sum_{\lambda \in \sigma^-(-\Delta+V)}|\lambda|^{\gamma}\leq
\frac{2^{\frac{d}{4}+1}}{(8\pi )^{\frac{d}{2}}}C_{HS}(\gamma)
\Big(\frac{ec}{2\gamma+d-2}\Big)^{\gamma+\frac{d}{2}-1}\frac{1}{[1-4\beta(c)]^{\frac{1}{2}}}\|V_-\|_{L^1}.
\end{equation*}
\end{theorem}
Similarly, from Theorem \ref{map02} we obtain
\begin{theorem}\label{mapw2}
Let $V$ be a Kato potential, and assume also $V_-\in L^2(\Real^d)$.
If $c> 0$ is such that $4\beta(c)<1$, then, for any $\gamma>2$,
\begin{equation*}\label{innw22}
\sum_{\lambda \in \sigma^-(-\Delta+V)}|\lambda|^{\gamma}\leq
\frac{2^{\frac{d}{4}}}{(8\pi )^{\frac{d}{2}}}C_{HS}(\gamma)
\Big(\frac{ec}{2\gamma+d-4}\Big)^{\gamma+\frac{d}{2}-2}\frac{1}{[1-4\beta(c)]^{\frac{1}{2}}}\|V_-\|_{L^2}^2.
\end{equation*}
\end{theorem}
We note that (\ref{al}) assures us that there always exists $c>0$
with $4\beta(c)<1$, so that Theorems \ref{mapw},\ref{mapw2} apply.

The dependence on $V_-$ in Theorems \ref{mapw},\ref{mapw2} is both
through its $L^1$-norm and through the quantity $\beta(c)$. The
quantity $\beta(c)$ can be written more explicitly by using the
integral representation of $(c-\Delta)^{-1}$,
$$((c-\Delta)^{-1}V_-)(x)=c^{\frac{d-2}{2}}\int_{\Real^d}G(c^{\frac{1}{2}}(x-y))V_-(y)dy,$$
where
$$G(x)=\frac{1}{(2\pi)^{\frac{d}{2}}}K_{\frac{d}{2}-1}(|x|)\frac{1}{|x|^{\frac{d}{2}-1}},$$
in which $K_{\frac{d}{2}-1}$ is the modified Bessel function of the
third kind (see, e.g., \cite{aron}). Thus
\begin{eqnarray}\label{bfor}\beta(c)&=&c^{\frac{d-2}{2}}\sup_{x\in \Real^d}
\int_{\Real^d}G(c^{\frac{1}{2}}(x-y))|V_-(y)|dy\nonumber\\
&=&\frac{1}{c}\sup_{x\in \Real^d}
\int_{\Real^d}G(x-y)|V_-(c^{-\frac{1}{2}}y)|dy.
\end{eqnarray}

We now introduce an apparently new norm on potentials, which is
natural in this context, in terms of which we can derive some useful
inequalities from Theorems \ref{mapw},\ref{mapw2}.
 For $\alpha> 0$, we say that a measurable
function $W:\Real^d\rightarrow \Real$ belongs to
$K^{\alpha}(\Real^d)$ if $\|W\|_{K^{\alpha}}<\infty$, where
\begin{eqnarray}\label{dew}\|W\|_{K^{\alpha}}&=&\sup_{c>0}c^{\alpha}\|(c-\Delta)^{-1}|W|\|_{L^{\infty}}\nonumber\\
&=&\sup_{x\in \Real^d,c>0}c^{\alpha-1}
\int_{\Real^d}G(x-y)|W(c^{-\frac{1}{2}}y)|dy.\end{eqnarray}
$K^{\alpha}(\Real^d)$ is a normed space with the above norm, and we
have $K^{\alpha}(\Real^d)\subset K(\Real^d)$ for all $\alpha>0$. By
the definition of the $K^\alpha$-norm and by (\ref{bfor}) we have,
when $V_-\in K^{\alpha}(\Real^d)$
\begin{equation}\label{betb}\beta(c)\leq
\|V_-\|_{K^{\alpha}}c^{-\alpha},\;\;\;\forall c>0.\end{equation}

To see that $K^\alpha(\Real^d)$ is a sufficiently large class of
functions, we note that
\begin{lemma}\label{ko}
If $d\geq 3$ and $p>\frac{d}{2}$ then $L^p(\Real^d)\subset
K^{\alpha}(\Real^d)$, where $\alpha=1-\frac{d}{2p}$, and we have,
for all $W\in L^p(\Real^d)$,
\begin{equation}\label{ui}\|W\|_{K^\alpha}\leq C_{d,p}\|
W\|_{L^p},
\end{equation}
where
\begin{equation}\label{dcdp}C_{d,p}=\Big(\int_{\Real^d}|G(x)|^{\frac{p}{p-1}}dx\Big)^{\frac{p-1}{p}}.
\end{equation}
\end{lemma}
\begin{proof}
Using H\"older's inequality we have
\begin{eqnarray*}\int_{\Real^d}G(x-y)|W(c^{-\frac{1}{2}}y)|dy
\leq
C_{d,p}c^{\frac{d}{2p}}\|W\|_{L^p}=C_{d,p}c^{1-\alpha}\|W\|_{L^p},
\end{eqnarray*}
which, using (\ref{dew}), implies (\ref{ui}). We note that the fact
that $C_{d,p}$ is
 finite follows from the condition $p>\frac{d}{2}$, which implies $\frac{p}{p-1}<\frac{d}{d-2}.$
\end{proof}

Another fact, which shows that $K^{\alpha}(\Real^d)$ contains
functions which are not in any $L^p(\Real^d)$ is
\begin{lemma}
If $W$ is measurable and
$$|W(x)|\leq  \frac{A}{|x|^{\eta}},\;\;\;\forall x\in\Real^d,$$
where $\eta\in (0,2)$, then $W\in K^{2-\eta}(\Real^d)$, and
$$\|W\|_{K^{2-\eta}}\leq \frac{1}{\pi^{\frac{d}{2}}}\frac{1}{2^{\eta+1}}\Gamma\Big(1-\frac{\eta}{2}
\Big)\Gamma\Big(\frac{d-\eta}{2} \Big)A$$
\end{lemma}
\begin{proof} We have
\begin{eqnarray*}&&\int_{\Real^d}G(x-y)|W(c^{-\frac{1}{2}}y)|dy \leq
Ac^{\frac{\eta}{2}}\int_{\Real^d}G(x-y)|y|^{-\eta}dy\\ &\leq&
Ac^{\frac{\eta}{2}}\int_{\Real^d}G(y)|y|^{-\eta}dy
=\frac{1}{\pi^{\frac{d}{2}}}\frac{1}{2^{\eta+1}}\Gamma\Big(1-\frac{\eta}{2}
\Big)\Gamma\Big(\frac{d-\eta}{2} \Big)Ac^{\eta-1},
\end{eqnarray*}
where the second inequality follows from the fact that both $G(x)$
and $|x|^{-\eta}$ are radially symmetric functions which are
decreasing in $|x|$, so that their convolution is maximized at the
origin.
\end{proof}

We now derive eigenvalue inequalities using the norms
$\|V_-\|_{K^{\alpha}}$. From Theorem \ref{mapw} and (\ref{betb}) we
have
\begin{equation}\label{innv11}
\sum_{\lambda \in \sigma^-(-\Delta+V)}|\lambda|^{\gamma}\leq
C_{HS}(\gamma)\frac{2^{\frac{d}{4}+1}}{(8\pi )^{\frac{d}{2}}}
\frac{\Big(\frac{ec}{2\gamma+d-2}\Big)^{\gamma+\frac{d}{2}-1}}{[1-4\|V_-\|_{K^{\alpha}}c^{-\alpha}]^{\frac{1}{2}}}\|V_-\|_{L^1}.
\end{equation}
We now wish to minimize the right-hand side of (\ref{innv11}) with
respect to $c$. We compute
$$\min_{
c>(4\|V_-\|_{K^{\alpha}})^{\frac{1}{\alpha}}}\frac{c^{\gamma+\frac{d}{2}-1}}{[1-4\|V_-\|_{K^{\alpha}}c^{-\alpha}]^{\frac{1}{2}}}
=\frac{2^\delta(2\delta+1)^{\delta+\frac{1}{2}}}{\delta^{\delta}}
 \|V_-\|_{K^{\alpha}}^{\delta},$$
 where
$$\delta=\frac{1}{\alpha}\Big(\gamma+\frac{d}{2}-1\Big).$$
Thus from (\ref{innv11}) we get
\begin{theorem}\label{cov2}
Let $V$ be a Kato potential, and assume also $V_-\in
L^1(\Real^d)\cap K^{\alpha}(\Real^d)$, where $\alpha>0$. Then, for
any $\gamma>2$,
\begin{equation*}
\sum_{\lambda \in \sigma^-(-\Delta+V)}|\lambda|^{\gamma}\leq
\kappa\|V_-\|_{L^1}\|V_-\|_{K^{\alpha}}^\delta,
\end{equation*}
where the constants are given by
$$\delta=\delta_{d,\alpha,\gamma}=\frac{1}{\alpha}\Big(\gamma+\frac{d}{2}-1\Big),$$
$$\kappa=\kappa_{d,\alpha,\gamma}=C_{HS}(\gamma)\frac{2^{\frac{d}{4}+1}}{(8\pi )^{\frac{d}{2}}}
\frac{2^{\delta}(2\delta+1)^{\delta+\frac{1}{2}}}{\delta^\delta}
\Big(\frac{e}{2\delta \alpha}\Big)^{\delta \alpha}.$$
\end{theorem}
Similarly, using Theorem \ref{mapw2} we obtain
\begin{theorem}\label{corv22}
Let $V$ be a Kato potential, and assume also $V_-\in
L^2(\Real^d)\cap K^\alpha(\Real^d)$, where $\alpha>0$. Then, for any
$\gamma>2$,
\begin{equation*}
\sum_{\lambda \in \sigma^-(-\Delta+V)}|\lambda|^{\gamma}\leq
\kappa\|V_-\|_{L^2}^2\|V_-\|_{K^\alpha}^{{\delta}},
\end{equation*}
where the constants are given by
$${\delta}={\delta}_{d,p,\gamma}=\frac{1}{\alpha}\Big(\gamma+\frac{d}{2}-2\Big),$$
$$\kappa={\kappa}_{d,\alpha,\gamma}=C_{HS}(\gamma)\frac{2^{\frac{d}{4}}}{(8\pi )^{\frac{d}{2}}}
\frac{2^\delta (2\delta+1)^{\delta+\frac{1}{2}}}{\delta^\delta}
 \Big(\frac{e}{2\delta
\alpha}\Big)^{\delta \alpha}.$$
\end{theorem}

We particularize to the case in which $d\geq 3$, $V_-\in
L^p(\Real^d)$, $p>\frac{d}{2}$. Using Lemma \ref{ko}, Theorem
\ref{cov2},\ref{corv22} imply
\begin{corollary}\label{cor2}
Assume $d\geq 3$. Let $V$ be a Kato potential, and assume also
$V_-\in L^1(\Real^d)\cap L^p(\Real^d)$, where $p>\frac{d}{2}$. Then,
for any $\gamma>2$,
\begin{equation}\label{res}
\sum_{\lambda \in \sigma^-(-\Delta+V)}|\lambda|^{\gamma}\leq
\kappa\|V_-\|_{L^1}\|V_-\|_{L^p}^\delta,
\end{equation}
where the constants are given by
\begin{equation}\label{ddel}\delta=\delta_{d,p,\gamma}=\frac{\gamma+\frac{d}{2}-1}{1-\frac{d}{2p}},\end{equation}
$$\kappa=\kappa_{d,p,\gamma}=C_{HS}(\gamma)(2C_{d,p})^{\delta}\frac{2^{\frac{d}{4}+1}}{(8\pi )^{\frac{d}{2}}}
\frac{(2\delta+1)^{\delta+\frac{1}{2}}}{\delta^\delta}
 \Big(\frac{e}{2\delta
(1-\frac{d}{2p})}\Big)^{\delta (1-\frac{d}{2p})},$$ with $C_{d,p}$
given by (\ref{dcdp}).
\end{corollary}

\begin{corollary}\label{cor22}
Assume $d\geq 3$. Let $V$ be a Kato potential, and assume also
$V_-\in L^2(\Real^d)\cap L^p(\Real^d)$, where $p>\frac{d}{2}$. Then,
for any $\gamma>2$,
\begin{equation}\label{res2}
\sum_{\lambda \in \sigma^-(-\Delta+V)}|\lambda|^{\gamma}\leq
\kappa\|V_-\|_{L^2}^2\|V_-\|_{L^p}^{{\delta}},
\end{equation}
where the constants are given by
$${\delta}={\delta}_{d,p,\gamma}=\frac{\gamma+\frac{d}{2}-2}{1-\frac{d}{2p}},$$
$$\kappa={\kappa}_{d,p,\gamma}=C_{HS}(\gamma)\Big(2C_{d,p} \Big)^{\delta}\frac{2^{\frac{d}{4}}}{(8\pi )^{\frac{d}{2}}}
\frac{(2\delta+1)^{\delta+\frac{1}{2}}}{\delta^\delta}
 \Big(\frac{e}{2\delta
(1-\frac{d}{2p})}\Big)^{\delta (1-\frac{d}{2p})},$$ with $C_{d,p}$
given by (\ref{dcdp}).
\end{corollary}
It is interesting to compare the inequalities given by Corollaries
\ref{cor2},\ref{cor22} with a different bound on the moments of
eigenvalues, given by the Lieb-Thirring inequalities
\cite{laptev,lieb}. These state that
\begin{equation}
\label{lieb} \sum_{\lambda \in
\sigma^-(-\Delta+V)}|\lambda|^{\gamma}\leq
C_{d,\gamma}\|V_-\|_{L^{\gamma+\frac{d}{2}}}^{\gamma+\frac{d}{2}},
\end{equation}
holds for any $\gamma\geq 0$ when $d\geq 3$, for any $\gamma>0$ when
$d=2$, and for any $\gamma\geq \frac{1}{2}$ when $d=1$.

Let us compare the bounds given by the inequalities when {\it{both}}
of them are valid. The following argument shows that our inequality
(\ref{res}) and the Lieb-Thirring inequality are independent, in the
sense that neither of them is stronger than the other: fixing
$\gamma>2$, $p>\frac{d}{2}$, if we take some potential $W\in
L^1(\Real^d)\cap L^{p}(\Real^d)\cap
L^{\gamma+\frac{d}{2}}(\Real^d)$, and define the family $V_{\mu}$
($\mu>0$) by
$$V_{\mu}(x)=\mu^{\frac{d}{\gamma+\frac{d}{2}}}W(\mu x)$$
then, for any $r>0$,
$$\|V_{\mu-}\|_{L^{r}}=\mu^{\frac{d}{\gamma+\frac{d}{2}}-\frac{d}{r}}\|W_-\|_{L^{r}},$$
hence
$$\|V_{\mu-}\|_{L^{\gamma+\frac{d}{2}}}^{\gamma+\frac{d}{2}}=\|W_-\|_{L^{\gamma+\frac{d}{2}}}^{\gamma+\frac{d}{2}},$$
\begin{eqnarray*}\|V_{\mu-}\|_{L^1}\|V_{\mu-}\|_{L^p}^{\delta}
=\mu^{-\frac{2d\delta}{(2\gamma+d)p}}
\|W_-\|_{L^1}\|W_-\|_{L^p}^{\delta},\end{eqnarray*} where $\delta$
is defined by (\ref{ddel}).
 Thus the
right-hand side of the inequality (\ref{res}) is arbitrarily small
for $\mu$ large and arbitrarily large for $\mu$ small, while the
right-hand side of (\ref{lieb}) does not depend on $\mu$, so that
our inequalities are sometimes weaker and sometimes stronger than
the Lieb-Thirring inequalities - depending on the potential $V$. In
particular (\ref{res}) is better than the bound given by the
Lieb-Thirring inequality when $\mu$ is large. A similar conclusion
holds with respect to the inequality (\ref{res2}) of Corollary
\ref{cor22}.

{\bf{Acknowledgement:}} We are grateful to M. Hansmann for his
critical reading of the manuscript and helpful comments.


\begin{thebibliography}{9}

\bibitem{aron} N. Aronszajn \& K.T. Smith, `Theory of Bessel
potentials. I.', Ann. Inst. Fourier {\bf{11}} (1961), 385-475.

\bibitem{krishna} M. Demuth \& M. Krishna, `Determining Spectra in
Quantum Theory', Birkh\"auser (Boston), 2005.

\bibitem{demuth} M. Demuth \& J.A. Van Casteren, Stochastic Spectral Theory for
Selfadjoint Feller Operators: A Functional Integration Approach,
Birkh\"auser (Basel), 2000.

\bibitem{laptev} A. Laptev \& T. Weidl,
Recent results on Lieb-Thirring inequalities, Journ\'ees \'Equations
aux d\'eriv\'ees partielles (2000) 1-14.

\bibitem{lieb} E.H. Lieb \& W. Thirring, Inequalities for the
moments of eigenvalues of the Schr\"odinger Hamiltonian and their
relation to Sobolev inequalities, Studies in Math. Phys., Essays in
honor of Valentine Bargmann, Princeton, 269-303 (1976).

\bibitem{rudin} W. Rudin, Real and Complex Analysis,
McGraw-Hill (New-York), 1987.

\bibitem{simonb} B. Simon, Trace Ideals and their Applications,
London Math. Soc. Lecture Notes, 1979.

\end{thebibliography}
\end{document}